\documentclass{amsproc}
\usepackage{amssymb}
\usepackage{graphicx}
\usepackage{amscd}
\usepackage{amsmath}
\theoremstyle{plain}

\newtheorem{corollary}{Corollary}

\newtheorem{definition}{Definition}

\newtheorem{proposition}{Proposition}
\newtheorem{remark}{Remark}

\newtheorem{theorem}{Theorem}
\numberwithin{equation}{section}

\def\cl{\operatorname{col}}

\begin{document}

\begin{center}

\noindent{\large{\bf Controllability of partial differential
equations on graphs}}
 \footnote{S.A.'s research is supported in part by the National Science Foundation,
grants ARC 0724860 and DMS-0648786 ; V.M's
 research is supported in part
by the University of Alaska Graduate Fellowship}

\bigskip

Sergei Avdonin and Victor Mikhaylov\\
Department of Mathematics and Statistics \\
     University of Alaska\\
     Fairbanks, AK 99775-6660, USA \\

\end{center}

\bigskip

\noindent{\footnotesize{\bf Abstract} We study the boundary
control problems for the wave, heat, and Schr\"odinger equations
on a finite graph. We suppose that the graph is a tree (i.e., it
does not contain cycles), and on each edge an equation is defined.
The control is acting through the Dirichlet condition applied to
all or all but one boundary vertices. The exact controllability in
$L_2$-classes of controls is proved and sharp estimates of the
time of controllability are obtained for the wave equation. The
null controllability for the heat equation and exact
controllability for the Schr\"odinger equation in arbitrary time
interval are obtained.

}

\vskip2mm \noindent{\footnotesize{\bf Classification Numbers}}:
\quad 81C05, 35R30, 35L05, 93B05, 49E15

\section{Introduction.}

Controllability problems for multi-link flexible structures or, in
other words, for the wave and beam equations on graphs were the
subject of extensive investigations of many mathematicians (see,
e.g. the review paper \cite{A} and references therein). Lagnese,
Leugering, and Schmidt in \cite{LLS,LLS2} used the method of
energy estimates together with the Hilbert Uniqueness Method to
show that the exact controllability can be achieved in optimal
time for tree-like graphs consisted of homogeneous strings, when
all but one exterior nodes are controlled. Independently Avdonin
and Ivanov \cite[Ch. VII]{AI} applied the method of moments and
the theory of vector-valued exponentials to study controllability
problems on graphs for the wave equation. The authors have proved
the exact controllability in the optimal time for the wave
equation on the star-shaped graph of non-homogeneous strings.
Belishev in \cite{BC,B} using the propagation of singularities
method obtained result on boundary controllability for a tree of
non-homogeneous strings with respect to the first component (the
shape) of the complete state.

The results on exact controllability fail as soon as cycles occur
within the network, even if all nodes (including the interior
nodes) are subjected to control \cite[Sec. VII.1]{AI}. However,
the spectral controllability may be retained for many graphs with
cycles (see \cite[Ch. VII]{AI}, \cite{DZB,LLS2} for details). In
\cite{LZ} for the tree of homogeneous vibrating strings, the
authors prove the exact controllability for some special class of
initial/final data. Many interesting results on spectral
controllability are obtained in \cite{DZB}.

In this paper we prove the result on the exact controllability for
the wave equation on a tree-like graph of non-homogeneous strings
for controls acting through Dirichlet conditions applied to all or
all but one boundary vertices. Our result generalizes the ones
from \cite{AI} and \cite{LLS2}. Using the controllability of the
wave equation and results from \cite{AI,M1,M2,R}, we prove also the null
controllability of the parabolic and exact controllability of the
Schr\"odinger equations on trees.

Controllability problems for partial differential equations on
graphs have many important applications. They are also related to inverse
problems on graphs \cite{AK,BC,BV} and to harmonic analysis \cite[Ch. VII]{AI}.
In this paper we use some known and prove several new results describing connections
between controllability of distributed parameter systems and properties of exponential families.

We do not consider some important problems closely related to the topic of the paper, such as
controllability of networks of beams and hybrid systems and refer a reader to the comprehensive papers
\cite{DN,LLS2,L}.

\section{Statement of the problems and main results.}

Let $\Omega$ be a finite connected compact graph without cycles (a
tree). The graph consists of edges $E=\{e_1,\ldots,e_N\}$
connected at the vertices $V=\{v_1\ldots,v_{N+1}\}$. Every edge
$e_j\in E$ is identified with an interval $(a_{2j-1},a_{2j})$ of
the real line. The edges are connected at the vertices $v_j$ which
can be considered as equivalence classes of the edge end points
$\{a_j\}$. The boundary $\Gamma=\{\gamma_1,\ldots,\gamma_m\}$ of
$\Omega$ is a set of vertices having multiplicity one (the
exterior nodes). We suppose that the graph is equipped with the
density
\begin{equation}
\label{density} \rho(x)\geqslant \operatorname{const}>0,\quad
x\in\Omega\backslash V, \quad \rho\in C^1(\bar e_j),\ j=1\ldots,N.
\end{equation}
All the results of this paper are valid also for piecewise continuously differentiable functions $\rho.$ Discontinuity of
$\rho$ or its derivative at the inner of an edge is equivalent to the additional inner vertex of
multiplicity two (see the compatibility conditions (\ref{Shr_cont}), (\ref{Kirh}) below).

 Since the graph under consideration is a tree, for
every $a,b\in\Omega,$ $a\not=b,$ there exist the unique path
$\pi[a,b]$ connecting these points. The density determines the
optical metric and the optical distance
\begin{eqnarray*}
d\,\sigma^2:=\rho(x)|d\,x|^2,\quad x\in\Omega\backslash V,\\
\sigma(a,b)=\int\limits_{\pi[a,\,b]}\sqrt{\rho(x)}|d\,x|, \quad
a,b\in\Omega,
\end{eqnarray*}
The optical diameter of the graph $\Omega$ is defined as
\begin{equation*}
\label{diam} d(\Omega)=\max_{a,\,b\in\Gamma}\sigma(a,b).
\end{equation*}
The graph $\Omega$ and the optical metric determine the {\it metric
graph} denoted by $\{\Omega,\rho\}$. For a rigorous definition of the metric
graph see, e.g. \cite{ks1,kuch,KuSt,AM1,graphbook}. The space of
real valued functions on the graph, square integrable with the
weight $\rho$ is denoted by $L_{2,\,\rho}(\Omega)$.

\subsection{Dirichlet spectral problem.}

Let $\partial w(a_j)$ denotes the derivative of $w$ at the vertex
$a_j$ taken along the corresponding edge in the direction toward
the vertex. We associate the following spectral problem to the
graph:
\begin{eqnarray}
-\frac{1}{\rho}\frac{d^2w}{dx^2}=\lambda w,\label{Shr_eqn}\\
w\in C(\Omega),\label{Shr_cont}\\
\sum_{a_j\in v}\partial w(a_j)=0\quad \text{for $v\in V\backslash\Gamma$}, \label{Kirh}\\
w=0\quad \text{on $\Gamma$}\label{Bound_Dir}.
\end{eqnarray}

It is well-known fact (see \cite{Be1,AM1,sol02}) that the problem
(\ref{Shr_eqn})--(\ref{Bound_Dir}) has a discrete spectrum of
eigenvalues $0<\lambda_1\leqslant\lambda_1\leqslant\lambda_2\ldots
$, $\lambda_k\to +\infty$ and corresponding eigenfunctions
$\phi_1,\phi_2,\ldots$ can be chosen so that
$\{\phi_k\}_{k=1}^\infty$ forms an orthonormal basis in
$\mathcal{H}:=L_{2,\,\rho}(\Omega)$:
\begin{eqnarray*}
(\phi_i,\phi_j)_\mathcal{H}=\int_\Omega
\phi_i(x)\phi_j(x)\rho(x)\,dx=\delta_{ij}\,.
\end{eqnarray*}
Set $\varkappa_k(\gamma)=\partial\phi_k(\gamma)$,
$\gamma\in\Gamma$. Let $\alpha_k$ be the $m$-dimensional column
vector defined as
$\alpha_k=\cl\left(\frac{\varkappa_k(\gamma)}{\sqrt{\lambda_k}}\right)_{\gamma\in\Gamma}.$
\begin{definition}
The set of pairs
\begin{equation}
\label{DSD} \left\{\lambda_k,\alpha_k\right\}_{k=1}^\infty
\end{equation}
is called the Dirichlet spectral data of the graph
$\{\Omega,\rho\}$.
\end{definition}

\subsection{Initial boundary value problems. Control from the whole boundary.}

We associate three dynamical systems, described correspondingly by
the wave, heat and Schr\"odinger equations, to the graph
$\{\Omega,\rho\}$. The first one has the form:
\begin{eqnarray}
\rho u_{tt}-u_{xx}=0\quad\text{in $\Omega\backslash V\times [0,T]$},\label{wv_eqn}\\
u|_{t=0}=u_t|_{t=0}=0,\label{cond_1}\\
u(\cdot,t)\quad \text{satisfies $(\ref{Shr_cont})$ and $(\ref{Kirh})$ for all $t\in[0,T]$},\label{cond_2}\\
u=f\quad\text{on $\Gamma\times [0,T]$}\label{cond_3}.
\end{eqnarray}
Here $T>0$, $f=f(\gamma,t)$, $\gamma\in\Gamma$, is the Dirichlet
boundary control which belongs to
$\mathcal{F}^T_\Gamma=L_2([0,T];\mathbb{R}^m)$. Inner product in
$\mathcal{F}^T_\Gamma$ is defined by
\begin{equation*}
(f,g)_{\mathcal{F}^T_\Gamma}=\sum_{i=1}^m\int_0^T
f(\gamma_i,t)g(\gamma_i,t)\,dt.
\end{equation*}

Let $D'(\Omega)$ be the set of distributions over the graph. We
introduce the space
\begin{equation*}
\mathcal{H}_{-1}=\left\{g\in D'(\Omega): \;g(x)=\sum_{k=1}^\infty
g_k\phi_k(x),\;
\left\{\frac{g_k}{\sqrt{\lambda_k}}\right\}_{k=1}^\infty\in
l_2\right\}.
\end{equation*}

The initial boundary value problem (\ref{wv_eqn})--(\ref{cond_3})
has a classical solution if $f\in C^2([0,T],\mathbb{R}^m)$. In our
case when $f\in \mathcal{F}^T_\Gamma$, the solution to
(\ref{wv_eqn})--(\ref{cond_3}) is understood in week
(distributional sense). It can be proved (see
\cite{AI,BV,DZB,LLS2}) that the solution $u^f$ satisfies the
inclusion
\begin{equation*} \label{reg}
u^f \in C([0,T];\mathcal{H} )  \cap
C^1([0,T]; \mathcal{H}_{-1}).
\end{equation*}
This means that $u(\cdot,t) \in \mathcal{H},\; u_t(\cdot,t) \in
\mathcal{H}_{-1}$ for all $t \in [0,T]$, and both functions are
continuous with respect to $t$ in corresponding norms. In other
words, the state of the dynamical system
(\ref{wv_eqn})--(\ref{cond_3}) $\left(u(\cdot,t),
u_t(\cdot,t)\right) $ is a point of $\mathcal{H} \times
\mathcal{H}_{-1},$ and the trajectory of the system is a
continuous curve in the state space $\mathcal{H} \times
\mathcal{H}_{-1}.$ This regularity result is sharp.

One of the main results of the present paper demonstrates the {\it
exact controllability} of the system
(\ref{wv_eqn})--(\ref{cond_3}).

\begin{theorem}
\label{main} For arbitrary state $\{a,b\}\in \mathcal{H}\times
\mathcal{H}_{-1}$, there exists such a control function
$f(\gamma,t)\in \mathcal{F}^T_\Gamma$ with $T=d(\Omega)$, that
solution of the initial boundary value problem
(\ref{wv_eqn})--(\ref{cond_3}) satisfies the equalities
$u^f(\cdot,T)=a$, $u^f_t(\cdot,T)=b$.
\end{theorem}

Another system we associate to the graph $\{\Omega,\rho\}$ is
\begin{eqnarray}
\rho u_t-
u_{xx}=0\quad\text{in $\Omega\backslash V\times [0,\tau]$},\label{par_eqn}\\
u|_{t=0}=a,\label{par_cond_1}\\
u(\cdot,t)\quad \text{satisfies $(\ref{Shr_cont})$ and $(\ref{Kirh})$ for all $t\in[0,\tau]$},\label{par_cond_2}\\
u=f\quad\text{on $\Gamma\times [0,\tau]$}\label{par_cond_3},
\end{eqnarray}
where $\tau>0$, $f\in \mathcal{F}^\tau_\Gamma$ and $a\in
\mathcal{H}_{-1}$.

It is known (see e.g. \cite{AI,Be3,LLS}), that the solution $u^f$
of the system (\ref{par_eqn})--(\ref{par_cond_3}) satisfies the
inclusion
\begin{equation*}
u^f \in C([0,\tau];\mathcal{H}_{-1}).
\end{equation*}

For the parabolic-type dynamical systems various types of
controllability are considered in the literature (see
\cite{AI,LLS2}). The following result demonstrates the {\it null
controllability} of the system
(\ref{par_eqn})--(\ref{par_cond_3}).

\begin{theorem}
\label{par_main} For arbitrary given state $a\in \mathcal{H}_{-1}$
and for arbitrary time interval $[0,\tau]$, $\tau>0$, there exists
a control $f\in \mathcal{F}^\tau_\Gamma$ such that the solution of
the initial boundary value problem
(\ref{par_eqn})--(\ref{par_cond_3}) satisfies the equality
$u^f(\cdot,\tau)=0$.
\end{theorem}

The Schr\"odinger equation can also be associated to the graph
$\{\Omega,\rho\}$:
\begin{eqnarray}
\rho u_t+
iu_{xx}=0\quad\text{in $\Omega\backslash V\times [0,\tau]$},\label{Schr_eqn}\\
u|_{t=0}=a,\label{Schr_cond_1}\\
u(\cdot,t)\quad \text{satisfies $(\ref{Shr_cont})$ and $(\ref{Kirh})$ for all $t\in[0,\tau]$},\label{Schr_cond_2}\\
u=f\quad\text{on $\Gamma\times [0,\tau]$}\label{Schr_cond_3},
\end{eqnarray}
where $f=f(\gamma,t)\in \mathcal{F}^T_\Gamma$,  $a\in
\mathcal{H}_{-1}$. It is known (see, e.g. \cite{ALP,TY}) that
solution $u^f(x,t)$ of (\ref{Schr_eqn})--(\ref{Schr_cond_3})
satisfies the inclusion
\begin{equation*}
u^f \in C([0,T];\mathcal{H}_{-1}).
\end{equation*}

For the dynamical system governed by the Schr\"odinger equation
(\ref{Schr_eqn})--(\ref{Schr_cond_3}) the following result on the
exact controllability holds. (Due to time reversibility, the exact
and null controllability are equivalent for the Schr\"odinger
equation.)
\begin{theorem}
\label{Schr_main} For arbitrary initial state $a\in \mathcal{H}_{-1}$
and for arbitrary time interval $[0,\tau]$, $\tau>0$, there exists
such a control $f\in \mathcal{F}^\tau_\Gamma$ that the solution to
the initial boundary value problem
(\ref{Schr_eqn})--(\ref{Schr_cond_3}) satisfies the equality
$u^f(\cdot,\tau)=0$.
\end{theorem}

\subsection{Initial boundary value problems. Control from a part of the boundary. }

In the case when the graph is controlled from the whole boundary
 but contains cycles,  the system
(\ref{wv_eqn})--(\ref{cond_3}) is not exactly controllable (see,
e.g. \cite[Sec. VII.1]{AI}). Similarly, if the graph is a tree,
but the system is not controlled at two or more  boundary points
(the Dirichlet condition $u=0$ is imposed there),  the Theorem
\ref{main} fails; the corresponding example (in the case of
homogeneous strings) is given in \cite[Sec. 6.3]{DZB} (see also
\cite{A}). Suppose that the graph is not controlled at one of the
boundary points, say $\gamma_1$. Then one can introduce the length
of the longest path from $\gamma_1$ to the rest of the boundary
$\Gamma_1=\Gamma \setminus \{\gamma_1\}$:
\begin{equation*}
d_1(\gamma_1,\Omega)=\max_{\gamma\in\Gamma_1}\tau(\gamma_1,\gamma).
\end{equation*}
The boundary conditions for the system
(\ref{wv_eqn})--(\ref{cond_2}) have the form:
\begin{equation}
\label{cond_4} u(\gamma_1,t)=0,\quad
u(\gamma_i,t)=f(\gamma_i,t),\quad i=2,\ldots,N,
\end{equation}
where $f\in \mathcal{F}^T_{\Gamma_1}=L_2([0,T];\mathbb{R}^{m-1})$.
In this situation the analog of Theorem \ref{main} holds true:

\begin{theorem}
\label{main_1} For arbitrary state $\{a,b\}\in \mathcal{H}\times
\mathcal{H}_{-1}$ there exists such a control function
$f(\gamma,t)\in \mathcal{F}^T_{\Gamma_1}$ with
$T=2d_1(\gamma_1,\Omega)$ that the solution of the initial
boundary value problem (\ref{wv_eqn})--(\ref{cond_2}),
(\ref{cond_4}) satisfies the equalities  $u^f(\cdot,T)=a$,
$u^f_t(\cdot,T)=b$.
\end{theorem}

For the parabolic and Schr\"odinger type systems
(\ref{par_eqn})--(\ref{par_cond_2}),
(\ref{Schr_eqn})--(\ref{Schr_cond_3}), we can also consider the
problem of the controllability from the part of the boundary,
i.e., we add the boundary conditions (\ref{cond_4}) to the
initial-value problem (\ref{par_eqn})--(\ref{par_cond_2}) and to the
problem (\ref{Schr_eqn})--(\ref{Schr_cond_2}). In this case one can
prove the analogs of Theorems \ref{par_main} and \ref{Schr_main}:
\begin{theorem}
\label{par_main_1} For arbitrary given state $a\in \mathcal{H}_{-1}$
and for arbitrary time interval $[0,\tau]$, $\tau>0$, there exists
such a control $f\in \mathcal{F}^\tau_{\Gamma_1}$ that the solution of
the initial boundary value problem
(\ref{par_eqn})--(\ref{par_cond_2}), (\ref{cond_4}) satisfies the
equality $u^f(\cdot,\tau)=0$.
\end{theorem}
\begin{theorem}
\label{Schr_main_1} For arbitrary initial state $a\in \mathcal{H}_{-1}$
and for arbitrary time interval $[0,\tau]$, $\tau>0$, there exists
such a control $f\in \mathcal{F}^\tau_{\Gamma_1}$ that the solution of
the initial boundary value problem
(\ref{Schr_eqn})--(\ref{Schr_cond_2}), $(\ref{cond_4})$ satisfies
the equality $u^f(\cdot,\tau)=0$.
\end{theorem}

\section{Auxiliary results.}

In \cite{BC}--\cite{BV} the following result concerning the controllability
with respect to the first component (the shape) of the complete
state $\{u,u_t\}$ of the dynamical system
(\ref{wv_eqn})--(\ref{cond_3}) has been proved:
\begin{theorem}
\label{Bel_control} Let $T=d(\Omega)/2$, then for arbitrary $a\in
\mathcal{H}$ there exists such a control $f(\gamma,t)\in
\mathcal{F}^T_\Gamma$ that the solution of the initial boundary
value problem (\ref{wv_eqn})--(\ref{cond_3}) satisfies the
equality $u^f(x,T)=a(x)$.
\end{theorem}
In other words, the system (\ref{wv_eqn})--(\ref{cond_3}) is
controllable with respect to the shape for the time equal to the
half optical diameter of the graph. Note that in general such a
control is not unique.

To prove Theorem {\ref{Bel_control}} the propagation of
singularities method has been used and the controllability was
reduced to solvability of the Volterra type equation.
It was supposed in \cite{BC}--\cite{BV} that  $\rho\in C^2$ on all edges, however, the
method works for $\rho\in C^1$ as well.
The same
technique can be applied to obtain the result on the
controllability of the system (\ref{wv_eqn})--(\ref{cond_3}) with
respect to the second component (the velocity) of the complete
state:
\begin{proposition}
\label{Teor_deriv} If $T=d(\Omega)/2$ then for arbitrary $b\in
\mathcal{H}_{-1}$, there exists such a control $f(\gamma,t)\in
\mathcal{F}^T_\Gamma$, that the solution of the initial boundary
value problem (\ref{wv_eqn})--(\ref{cond_3}) satisfies the
equality $u_t^{f}(x,T)=b(x)$.
\end{proposition}

In the following two propositions we consider the case of boundary
condition $(\ref{cond_4})$ for the system
(\ref{wv_eqn})--(\ref{cond_2}). The proof of the first proposition
can be extracted from the proof of Theorem $\ref{Bel_control}$
 \cite[Sec. 2]{BV}. Let us introduce the `optical center' of
the graph $\Omega$, i.e., such a point $\xi\in\Omega$, that
$\max\limits_{\gamma\in\Gamma}\tau(\xi,\gamma)=d(\Omega)/2=T$.
Since $\Omega$ is a tree, there can be only one optical center.
Suppose that the final state $a(x)$ is supported in such a subtree
$\Omega_1\subset\Omega$ that $\xi\notin\Omega_1$. As it was shown
in \cite{BC}--\cite{BV}, to solve the control problem one need to use
controls supported on the part of the boundary of the graph
$\Omega$ which is the boundary of $\Omega_1$. In other words, it
is possible to construct such a control $f\in
\mathcal{F}^T_\Gamma$ that $u^f(T,x)=a(x)$ and $f(\gamma,t)=0$ for
$\gamma\notin\Omega_1$. The authors offers the explicit procedure
of the construction of such a control. If instead of the `optical
center' of the graph we take a boundary point $\gamma_1$ where the
homogeneous Dirichlet condition $u(\gamma_1,t)=0$ is imposed, we
come to the following statements:

\begin{proposition}
\label{rem_2} If $T=d_1(\gamma_1,\Omega)$, then for arbitrary
$a\in \mathcal{H}$, there exists such a control $f\in
\mathcal{F}^T_{\Gamma_1}$, that the solution of the boundary value
problem (\ref{wv_eqn})--(\ref{cond_2}), (\ref{cond_4}) satisfies
the equality $u^f(x,T)=a(x)$.
\end{proposition}
The same result holds true for the controllability with respect to
the velocity:
\begin{proposition}
\label{rem_3} If $T=d_1(\gamma_1,\Omega)$ then for arbitrary $b\in
\mathcal{H}_{-1}$, there exists such a control $f\in
\mathcal{F}^T_{\Gamma_1}$ that the solution of the boundary value
problem (\ref{wv_eqn})--(\ref{cond_2}), (\ref{cond_4}) satisfies
the equality $u^f_t(x,T)=b(x)$.
\end{proposition}

\section{Proof of Theorem $\ref{main}$.}

We begin with the reducing the problem of controllability of the
dynamical system (\ref{wv_eqn})--(\ref{cond_3}) to the moment
problem in $\mathcal{F}^T_\Gamma$. Solving the initial boundary
value problem (\ref{wv_eqn})--(\ref{cond_3}) by the Fourier method
and looking for the solution in the form
\begin{equation}
\label{Four_repr} u^f(x,t)=\sum_{k=1}^\infty c_k^f(t)\phi_k(x),
\end{equation}
we get the expression for the coefficients:
\begin{equation*}
\label{coeff}
c_k^f(t)=\sum_{\gamma\in\Gamma}\frac{\varkappa_k(\gamma)}{\sqrt{\lambda_k}}\int_0^t\sin{\sqrt{\lambda_k}(t-s)}
f(\gamma,s)\,ds.
\end{equation*}
Suppose that we are given the final state
$\{a,b\}\in\mathcal{H}\times \mathcal{H}_{-1}$ at $t=T$, where the
functions $a(x)$, $b(x)$ have the expansions
\begin{eqnarray*}
a(x)=\sum_{k=1}^\infty a_k\phi_k(x),\quad b(x)=\sum_{k=1}^\infty
b_k\phi_k(x),
\end{eqnarray*}
for some $\{a_k\}_{k=1}^\infty\in l_2$ and
$\{\frac{b_k}{\sqrt{\lambda_k}}\}_{k=1}^\infty\in l_2$. Then for
an unknown control $f\in\mathcal{F}^T_\Gamma,$ the following moment
equalities should hold at time $t=T$:
\begin{eqnarray}
a_k=c_k^f(T)=\sum_{\gamma\in\Gamma}\frac{\varkappa_k(\gamma)}{\sqrt{\lambda_k}}
\int_0^T\sin{\sqrt{\lambda_k}(T-s)}f(\gamma,s)\,ds,\quad k\in \mathbb{N},\label{moment_1}\\
\frac{b_k}{\sqrt{\lambda_k}}=\frac{\dot
c_k^f(T)}{\sqrt{\lambda_k}}=\sum_{\gamma\in\Gamma}\frac{\varkappa_k(\gamma)}{\sqrt{\lambda_k}}\int_0^T\cos{\sqrt{\lambda_k}(T-s)}
f(\gamma,s)\,ds,\quad k\in \mathbb{N}\label{moment_2}.
\end{eqnarray}
Using Euler formulas for exponentials, we rewrite
$(\ref{moment_1})$, $(\ref{moment_2})$ as
\begin{eqnarray}
\frac{b_k}{\sqrt{\lambda_k}}\pm
ia_k=\sum_{\gamma\in\Gamma}\frac{\varkappa_k(\gamma)}{\sqrt{\lambda_k}}\int_0^T
e^{\pm i\sqrt{\lambda_k}(T-s)}f(\gamma,s)\,ds,\quad k\in
\mathbb{N}.\label{mom_1}
\end{eqnarray}
\begin{definition}
We call the moment problem $(\ref{mom_1})$ solvable (and
$f(\gamma,t)$ a solution of the moment problem) in the time
interval $[0,T]$ if, for arbitrary $\{a_k\}_{k=1}^\infty$,
$\{\frac{b_k}{\sqrt{\lambda_k}}\}_{k=1}^\infty\in l_2$, there
exist such a function $f\in \mathcal{F}^T_\Gamma$ that equalities
$(\ref{mom_1})$ hold.
\end{definition}
We emphasize that the solvability of the moment problem
$(\ref{mom_1})$ in the time interval $[0,T]$ for some $T>0$ is
equivalent to the controllability of the dynamical system
(\ref{wv_eqn})--(\ref{cond_3}) in the sense of Theorem \ref{main}
in the same time interval. This is a basic statement of the method
of moments (see, e.g. \cite[Ch. III]{AI}, \cite{R}).

We need a couple of definitions concerning vector families in
arbitrary Hilbert space.
\begin{definition}
The family $\{\xi_k\}_{k=1}^\infty$ in a Hilbert space $H$ is
called a Riesz basis, if it is an image of an orthonormal basis
under the action of some linear isomorphism.
\end{definition}

\begin{definition}
The family $\{\xi_k\}_{k=1}^\infty$ in a Hilbert space $H$ is
called an $\mathcal{L}$-basis, if it is a Riesz basis in the
closure of the linear span of the family.
\end{definition}

The result on the controllability formulated in Theorem
\ref{Bel_control} implies the solvability of the moment problem
$(\ref{moment_1})$ for $T=d(\Omega)/2$ for every
$\{a_k\}_{k=1}^\infty$. The controllability result formulated in
Proposition \ref{Teor_deriv} implies the solvability of the moment
problem $(\ref{moment_2})$ for $T=d(\Omega)/2$ for every
$\{\frac{b_k}{\sqrt{\lambda_k}}\}_{k=1}^\infty\in l_2$. Our goal
is to show that the solvability of the moment problems
$(\ref{moment_1})$ and $(\ref{moment_2})$ for $T=d(\Omega)/2$
implies the solvability of the moment problem (\ref{mom_1}) for
$T=d(\Omega)$. Let us put $T_*=d(\Omega)/2$ and introduce the
families of vector valued functions
\begin{eqnarray*}
S_k(t)=\alpha_k\sin{\sqrt{\lambda_k}t},\quad
C_k(t)=\alpha_k\cos{\sqrt{\lambda_k}t},\quad k\in \mathbb{N}.
\end{eqnarray*}

According to Theorem III.3.3 of \cite{AI} the solvability of the
moment problems $(\ref{moment_1})$ and $(\ref{moment_2})$ means
that the families $\left\{S_k(t) \right\}_{k=1}^\infty$ and
$\left\{C_k(t) \right\}_{k=1}^\infty $ form  $\mathcal{L}$-bases
in $L_2([0,T_*];\mathbb{R}^m)$.

Let us introduce subspaces of $L_2([0,T_*]; \mathbb{R}^m)$:
\begin{eqnarray*}
\Xi_o=\bigvee\{S_k(t)\}_{k=1}^\infty,\quad
\Xi_e=\bigvee\{C_k(t)\}_{k=1}^\infty,
\end{eqnarray*}
where $\bigvee$ denotes the closure of the linear span of a family. We extend
the functions from $\Xi_o$ to the interval $[-T_*,0)$ in the odd
way:
\begin{eqnarray*}
\widetilde\varphi(t)=\left\{\begin{array}l \varphi(t),\quad t\geqslant 0,\\
-\varphi(-t),\quad t< 0,\end{array}\right.,\quad -T_*\leqslant
t\leqslant T_*,\quad \varphi\in\Xi_o\,,
\end{eqnarray*}
and the functions from $\Xi_e$ --- in the even way:
\begin{eqnarray*}
\widetilde\varphi(t)=\left\{\begin{array}l \varphi(t),\quad t\geqslant 0,\\
\varphi(-t),\quad t< 0,\end{array}\right.\quad -T_*\leqslant
t\leqslant T_*,\quad \varphi\in\Xi_e\,.
\end{eqnarray*}
Let us denote the spaces of extended functions by
$\widetilde\Xi_o$ and $\widetilde\Xi_e$ and notice that the
extended families $\{\widetilde{S}_k(t)\}_{k=1}^\infty$ and
$\{\widetilde{C}_k(t)\}_{k=1}^\infty$ are Riesz bases in
$\widetilde\Xi_o$ and $\widetilde\Xi_e$ correspondingly. The
orthogonality of the spaces $\widetilde\Xi_o$ and
$\widetilde\Xi_e$ implies that the union
\begin{equation*}
\Bigl\{\widetilde{C}_k(t)\Bigr\}_{k=1}^\infty\bigcup
\Bigl\{\widetilde{S}_k(t)\Bigr\}_{k=1}^\infty\,,
\end{equation*}
forms a Riesz basis in $\widetilde\Xi_o\bigoplus
\widetilde\Xi_e\subset L_2([-T_*,T_*];
\mathbb{R}^m)$. 
Introducing functions
\begin{eqnarray}
\label{expon} E_{\pm k}(t)=C_k(t)\pm iS_k(t)=\alpha_ke^{\pm
i\sqrt{\lambda_k}t},\quad k\in \mathbb{N},
\end{eqnarray}
we see that the set $\{E_{\pm k}\}_{k\in \mathbb{N}}$ forms an
$\mathcal{L}$-basis in $L_2([-T_*,T_*]; \mathbb{C}^m)$. Shifting
the argument, we come to the conclusion that the  same family
forms an $\mathcal{L}$-basis in $L_2([0,2T_*]; \mathbb{C}^m)$.
Then according to Theorem III.3.3 of \cite{AI}, the moment problem
$(\ref{mom_1})$ is solvable for time $T=2T_*=d(\Omega)$. As we
have already noticed, this implies the exact controllability of
(\ref{wv_eqn})--(\ref{cond_3}) in the time interval
$[0,d(\Omega)]$. Theorem $\ref{main}$ is proved.

The proof of Theorem $\ref{main_1}$ is analogous to the previous one.
We set $\alpha'_k$ to be the $(m-1)$-dimensional column vector
defined as
\begin{equation}
\label{DSD_1}
\alpha_k'=\cl\left(\frac{\varkappa(\gamma)}{\sqrt{\lambda_k}}\right)_{\gamma\in\Gamma_1}.
\end{equation}
There naturally arise the families of vector functions in
$L_2([0,T_1]; \mathbb{R}^{m-1})$ with $T_1=d_1(\gamma_1,\Omega)$:
\begin{eqnarray*}
S^1_k(t)=\alpha_k'\sin{\sqrt{\lambda_k}t},\quad
C^1_k(t)=\alpha_k'\cos{\sqrt{\lambda_k}t},\quad k\in \mathbb{N}.
\end{eqnarray*}
One should perform the same procedure (using Propositions
$\ref{rem_2}$ and $\ref{rem_3}$ instead of Theorem
\ref{Bel_control} and Proposition \ref{Teor_deriv}) as in the
proof of Theorem $\ref{main}$, construct the family of vector
exponentials
\begin{equation}
\label{expon_1} \left\{E^1_{\pm k}\right\}_{k\in N},\quad E^1_{\pm
k}(t)=\alpha_k'e^{\pm i\sqrt{\lambda_k}t},\quad t\in
(0,2T_1),\quad k\in \mathbb{N},
\end{equation}
and use the connection between controllability and vector
exponentials (\cite{AI}, Theorem III.3.3).

In the proofs of Theorems \ref{main} and \ref{main_1} we have got
important results which are of independent interest in Function
Theory.

\begin{proposition} \label{wb}The family
$\left\{E_{\pm k}\right\}_{k\in \mathbb{N}}$ (see $(\ref{expon})$)
constructed using the Dirichlet spectral data $(\ref{DSD})$ is an
$\mathcal{L}$-basis in $L_2([0,d(\Omega)]; \mathbb{C}^m)$.
\end{proposition}
Suppose that we pick arbitrary boundary point of the graph (we
keep the notation $\gamma_1$ for it), then we get
\begin{proposition}
The family $\left\{E^1_{\pm k}\right\}_{k\in \mathbb{N}}$ (see
$(\ref{expon_1})$) constructed using the Dirichlet spectral data
$(\ref{DSD})$, $(\ref{DSD_1})$ is an $\mathcal{L}$-basis in
$L_2([0,2T_1]; \mathbb{C}^{m-1})$ for $T_1=d_1(\gamma_1,\Omega)$.
\end{proposition}
It seems to be very difficult to obtain these results without
using the control theoretic approach.

\section{Proof of Theorem $\ref{par_main}$.}

Looking for the solution of (\ref{par_eqn})--(\ref{par_cond_3}) in
the form $(\ref{Four_repr})$ for the fixed initial state $a\in
\mathcal{H}_{-1}$ with the expansion
\begin{equation}
\label{Par_init_expan}
a(x)=\sum\limits_{k=1}^\infty a_k\phi_k(x),
\end{equation}
we come to the following formulas for the coefficients:
\begin{equation*}
\label{par_coeff}
c_k^f(t)=a_ke^{-\lambda_kt}+\sum_{\gamma\in\Gamma}\varkappa_k(\gamma)
\int_0^te^{-\lambda_k(t-s)} f(\gamma,s)\,ds.
\end{equation*}
Solving the control problem associated with
(\ref{par_eqn})--(\ref{par_cond_3}) in the time interval
$[0,\tau]$, we need the equation $c_k^f(\tau)=0$, $k\in
\mathbb{N}$ to be satisfied. This leads to the following moment
problem
\begin{equation}
\label{par_mom}
0=\frac{a_k}{\sqrt{\lambda_k}}e^{-\lambda_k\tau}+\sum_{\gamma\in\Gamma}\frac{\varkappa_k(\gamma)}{\sqrt{\lambda_k}}\int_0^\tau
e^{-\lambda_k(\tau-s)} f(\gamma,s)\,ds,\quad k\in \mathbb{N}.
\end{equation}

\begin{definition}
\label{Def_moment_sol} The moment problem $(\ref{par_mom})$ is
solvable in the time interval $[0,\tau]$ for some $\tau>0$ if and
only if, for arbitrary
$\left\{\frac{a_k}{\sqrt{\lambda_k}}\right\}_{k=1}^\infty\in l_2,$
there exists $f\in \mathcal{F}^\tau_\Gamma$ such that equalities
$(\ref{par_mom})$ hold.
\end{definition}
Note that solvability of the moment problem (\ref{par_mom}) is
equivalent to the null controllability of the dynamical system
(\ref{par_eqn})--(\ref{par_cond_3}).
\begin{definition}
The family $\{\xi_k\}_{k=1}^\infty$ in a Hilbert space $H$ is
called minimal if, for every $k\in \mathbb{N}$, element $\xi_k$
does not belong to the closure of the linear span of the remaining
elements.
\end{definition}
Another equivalent characteristic of the minimal family
$\{\xi_k\}_{k=1}^\infty$ in a Hilbert space $H$ with the scalar
product $<\cdot,\cdot>$ is the existence of the bi-orthogonal
family $\{\xi_k'\}_{k=1}^\infty\subset H$ such that
\begin{equation*}
<\xi_k,\xi_n'>=\delta_{k,n}, \quad k,n\in \mathbb{N}.
\end{equation*}

It is well known, that if a vector family is an
$\mathcal{L}$-basis in $H$, it is minimal in $H$.

Proposition \ref{wb} states that the `hyperbolic' family $\{E_{\pm k}\}_{k\in \mathbb{N}}$
defined by $(\ref{expon})$ forms an
$\mathcal{L}$-basis in  $L_2([0,d(\Omega)]; \mathbb{C}^m)$. Let
us denote by $\{E_{\pm k}'\}_{k\in \mathbb{N}}$ the family
bi-orthogonal to $\{E_{\pm k}\}_{k\in \mathbb{N}}$. There are connections between the
`hyperbolic' family $(\ref{expon})$ and the `parabolic' one,
\begin{equation} \label{par_expon}
\left\{Q_k\right\}_{k=1}^\infty,\quad Q_ k(t)=\alpha_ke^{-
\lambda_kt},\quad k\in \mathbb{N},
\end{equation}
first established by D.L. Russell \cite{R}. We use his result in a
slightly more general form, formulated in Theorem II.5.20 of
\cite{AI}, from which it follows that the `parabolic' family
$\left\{Q_k\right\}_{k=1}^\infty$ is minimal in $L_2([0,\tau],
\mathbb{C}^m)$ for every $\tau>0$ and for the members of the
`parabolic' bi-orthogonal family $\{Q_k'\}_{k=1}^\infty$ the
following estimates hold:
\begin{equation}
\label{Biort_est} \|Q_k'\|_{L_2([0,\tau],
\mathbb{C}^m)}\leqslant C(\tau)
\|E'_k\|_{L_2([0,d(\Omega)],\mathbb{C}^m)}e^{\beta\sqrt{|\lambda_n|}},\quad k\in
\mathbb{N},
\end{equation}
with positive constants $C(\tau)$ and $\beta$.

To prove Theorem \ref{par_main}, one needs to show the solvability
of the moment problem $(\ref{par_mom})$ which  can be rewritten
as
\begin{equation*} \label{par_mom_1}
-\frac{a_k}{\sqrt{\lambda_k}}e^{-\lambda_k\tau} = \sum_{\gamma\in\Gamma}\frac{\varkappa_k(\gamma)}{\sqrt{\lambda_k}}\int_0^\tau
e^{-\lambda_kt} f(\gamma,\tau-t)\,dt,\quad k\in \mathbb{N},
\end{equation*}
or, shortly, as
\begin{equation} \label{par_mom_2} -\frac{a_k}{\sqrt{\lambda_k}}e^{-\lambda_k\tau}=\left(
Q_k,f^\tau \right)_{\mathcal{F}^\tau_\Gamma}, \quad k\in \mathbb{N},\end{equation}
where $f^\tau(\gamma,t)=f(\gamma,\tau-t).$
One can check that a formal solution of (\ref{par_mom_2}) has the form
\begin{equation}
\label{par_contr_sol} f^\tau(\gamma,t)=-\sum\limits_{k=1}^\infty
a_ke^{-\lambda_k\tau}Q_k'(t).
\end{equation}
Estimates $(\ref{Biort_est})$ imply that $f^\tau(\gamma,t)$ defined by
$(\ref{par_contr_sol})$, belongs to $\mathcal{F}^\tau_\Gamma$ and
therefore, the moment problem $(\ref{par_mom})$ is solvable. This completes
the proof of Theorem \ref{par_main}.

The proof of Theorem $\ref{par_main_1}$ is similar. The
corresponding family of exponentials that arise while reducing the
control problem to the moment problem has the form:
\begin{equation}
\label{par_expon_1} \left\{Q_k^1\right\}_{k=1}^\infty,\quad Q_
k^1(t)=\alpha_k'e^{- \lambda_kt}\quad k\in \mathbb{N}.
\end{equation}

We conclude this section with results about families of vector
exponentials that naturally appeared in the proofs.
\begin{proposition}
\label{Cor_par_1} The family $\left\{Q_k\right\}_{k=1}^{\infty}$
(see $(\ref{par_expon})$) constructed using the Dirichlet spectral
data $(\ref{DSD})$ is minimal in $L_2([0,T]; \mathbb{C}^m)$ for
any $T>0$.
\end{proposition}
Suppose that we pick arbitrary boundary point of the graph (we
keep the notation $\gamma_1$ for it), then the following statement is true.
\begin{proposition}
\label{Cor_par_2} The family $\left\{Q^1_k\right\}_{k=1}^{\infty}$
(see $(\ref{par_expon_1})$) constructed using the Dirichlet
spectral data $(\ref{DSD})$, $(\ref{DSD_1})$ is minimal in
$L_2([0,T]; \mathbb{C}^{m-1})$ for any $T>0$.
\end{proposition}
We emphasize that an independent proof of Propositions
\ref{Cor_par_1}, \ref{Cor_par_2} without using the control
theoretic approach would be very difficult.

\section{Proof of Theorem $\ref{Schr_main}$.}

To prove Theorem \ref{Schr_main} we use the scheme proposed in
\cite{M1}. We reformulate the initial boundary value problems
(\ref{wv_eqn})--(\ref{cond_3}) and
(\ref{Schr_eqn})--(\ref{Schr_cond_3}) in the operator form.
Results concerning the dependance of solutions to systems dual to
(\ref{wv_eqn})--(\ref{cond_3}),
(\ref{Schr_eqn})--(\ref{Schr_cond_3}) on the initial data, allow
us to use the Theorem $3.1$ of \cite{M1} that derives the exact
controllability of the first-order system
(\ref{Schr_eqn})--(\ref{Schr_cond_3}) in the arbitrary time
interval from the exact controllability of the second-order system
(\ref{wv_eqn})--(\ref{cond_3}) in some time interval.

Let us introduce the operator $A=-\frac{1}{\rho}\frac{d^2}{dx^2}$
in $H_0:=\mathcal{H}=L_{2,\ \rho}(\Omega)$. If the density $\rho$
satisfies (\ref{density}), the operator $A$ is self-adjoint,
positive definite and boundedly invertible with the domain
\begin{equation*}
D(A)=\left\{a\in H_0,\, a|_{e_i}\in H^2(e_i),\, \text{$a$
satisfies (\ref{Shr_cont}), (\ref{Kirh}))},\, a|_\Gamma=0\right\}.
\end{equation*}
This operator defines the scale $H_p$, $p\in \mathbb{Z}$, of
Hilbert spaces. For $p>0$, integer, $H_p=D(A^{\frac{p}{2}})$ with
the norm $\|x\|_p=\|A^{\frac{p}{2}}x\|$, $H_{-p}$ is dual to $H_p$
with respect to the scalar product in $H_0$. Another
characterization of $H_{-p}(\Omega)$ is that it is the completion
of $H_0$ with respect to the norm
$\|x\|_{-p}=\|A^{-\frac{p}{2}}x\|$. 
By $A'$ we denote the operator dual to $A$: it is the extension of
$A$ to $H_{-2}$ with the domain $H_0$. Let $Y=\mathbb{R}^m$ and
$C: H_2\mapsto Y$ be defined by:
\begin{equation*}
Ca=\cl\left(\partial a(\gamma)\right)_{\gamma\in\Gamma}.
\end{equation*}
Let the operator $B:Y\mapsto H_{-2}$ be dual to $C$. In this
notations we can rewrite the dynamical system
(\ref{Schr_eqn})--(\ref{Schr_cond_3}) as
\begin{equation}
\label{Schr_op} u_t(t)-iA'u(t)=Bf(t),\quad u(0)=a\in H_0.
\end{equation}
The dual observation system with output function $y$ is defined by
\begin{equation}
\label{Schr_op_dual} u_t(t)-iAu(t)=0,\quad u(0)=u_0\in H_0,\quad
y(t)=Cu(t).
\end{equation}
The smoothness of the solution of (\ref{Schr_op_dual}) (see
\cite{ALP} for the case of one interval) guarantees that for the
{\it observation} operator $\mathcal{C}_s: u_0\mapsto y(t)$ the
following estimate holds:
\begin{equation}
\label{Schr_admiss} \|\mathcal{C}_su_0\|_{\mathcal{F}^T}\leqslant
K_T\|u_0\|_{H_0}, \quad u_0\in H_2
\end{equation}
with $K_T>0$. 

System (\ref{wv_eqn})--(\ref{cond_3}) can be rewritten as
\begin{equation}
\label{Wave_op} u_{tt}(t)+A'u(t)=Bf(t),\quad u(0)=0,\, u_t(0)=0.
\end{equation}
The dual observation system with the output function $z$ has the form
\begin{equation*}
\label{Wave_op_dual} u_{tt}(t)+Au(t)=0,\quad u(0)=u_0\in H_1,\,
u_t(0)=u_1\in H_0,\quad z(t)=Cu(t).
\end{equation*}
The {\it observation} operator
$\mathcal{C}_w:\{u_0,u_1\}\mapsto z(t),$ satisfies the estimate:
\begin{equation}
\label{wave_admiss}
\|\mathcal{C}_w\{u_0,u_1\}\|_{\mathcal{F}^T}\leqslant
K^1_T(\|u_0\|_{H_1}+\|u_1\|_{H_0})
\end{equation}
with $K_T^1>0$ (see
\cite{LLT}).
Now we can use Theorem $3.1$ of \cite{M1}, which says that if the
dynamical system (\ref{Wave_op}) is exactly controllable in some
time interval (in our case it is controllable in the time interval
$(0,d(\Omega)$), then the system (\ref{Schr_op}) is exactly
controllable in any time interval, provided observation operators
satisfy inequalities (\ref{Schr_admiss}), (\ref{wave_admiss}).
This completes the proof of Theorem \ref{Schr_main}.

\begin{remark}
The proof of Theorem \ref{Schr_main_1} is similar, one should
refer to Theorem \ref{main_1} for the controllability of the
corresponding second order dynamical system.
\end{remark}

Looking for the solution of (\ref{Schr_eqn})--(\ref{Schr_cond_3})
in the form $(\ref{Four_repr})$ for the fixed initial state $a\in
\mathcal{H}_{-1}$ with the expansion (\ref{Par_init_expan}),
we come to the following formulas for the coefficients:
\begin{equation*}
\label{Schr_coeff}
c_k^f(t)=a_ke^{i\lambda_kt}+\sum_{\gamma\in\Gamma}
{\varkappa_k}(\gamma)\int_0^te^{i\lambda_k(t-s)} f(\gamma,s)\,ds.
\end{equation*}
Solving the control problem associated with
(\ref{Schr_eqn})--(\ref{Schr_cond_3}) in the time interval
$[0,\tau]$, we obtain the following moment problem
\begin{equation}
\label{Schr_mom}
0=\frac{a_k}{\sqrt{\lambda_k}}e^{i\lambda_k\tau}+\sum_{\gamma\in\Gamma}
\frac{{\varkappa_k(\gamma)}}{\sqrt{\lambda_k}}\int_0^\tau
e^{i\lambda_k(\tau-s)} f(\gamma,s)\,ds,\quad k\in \mathbb{N}.
\end{equation}
Theorem \ref{Schr_main} implies that the moment problem
(\ref{Schr_mom}) is solvable for any $\tau>0$. Using Theorem
III.3.3 of \cite{AI} we deduce the result about family of vector
valued exponentials that appeared in the moment problem
(\ref{Schr_mom}).
\begin{corollary}
The family
\begin{equation*}
\left\{D_k\right\}_{k=1}^\infty,\quad D_ k(t)= {\alpha_k}e^{i
\lambda_kt},\quad k\in \mathbb{N},
\end{equation*}
constructed using the Dirichlet spectral data (\ref{DSD}) is an
$\mathcal{L}$-basis in $L_2([0,\tau]; \mathbb{C}^m)$ for any
$\tau>0$.
\end{corollary}
Picking arbitrary boundary point of the graph (we keep the
notation $\gamma_1$ for it) and  using Theorem
\ref{Schr_main_1}, we get
\begin{corollary}
The family
\begin{equation*}
\left\{D^1_k\right\}_{k=1}^\infty,\quad D^1_ k(t)= \alpha_k'e^{i
\lambda_kt},
\end{equation*}
constructed using the Dirichlet spectral data (\ref{DSD}),
(\ref{DSD_1}) is an $\mathcal{L}$-basis in $L_2([0,\tau];
\mathbb{C}^{m-1})$ for any $\tau>0$.
\end{corollary}

\section{Acknowledgments}

The authors are very grateful to the referees for valuable
remarks and suggestions and to M.~Belishev and P.~Kuchment for fruitful discussions.

\end{document}